\documentclass[titlepage,openany,rotate]{amsart}

\usepackage{amsmath}
\usepackage{amssymb, amsfonts,amscd}

\theoremstyle{plain}
\newtheorem{Prop}{Proposition}[section]
\newtheorem{Thm}[Prop]{Theorem}
\newtheorem{Cor}[Prop]{Corollary}
\newtheorem{Lem}[Prop]{Lemma}

\theoremstyle{definition}
\newtheorem{Def}[Prop]{Definition}

\theoremstyle{remark}
\newtheorem{Rem}[Prop]{Remark}
\newtheorem{Problem}[Prop]{\bf Problem}
\newtheorem{Example}[Prop]{\bf Example}
\newtheorem{Question}[Prop]{\bf Question}

\def\dim{\mathop{\roman{dim}}}
\def\int{\mathop{\roman{int}}}

\def\id{\mathop{\roman{id}}}

\def\1{^{-1}}

\def\NN{{\bf N}}

\def\dim{\text{dim}}
\def\id{\text{id}}

\def\Int{\text{Int}}
\def\Cone{\text{Cone}}
\def\Tel{\mbox{Tel}}
\def\lqs{\le}

\def\CCC{{\mathcal C}}
\def\FFF{{\mathcal F}}
\def\GGG{{\mathcal G}}

\def\dokaz{{\bf Proof. }}
\def\edokaz{\hfill $\blacksquare$}
\newcommand\dokazp[1]{\vspace{\baselineskip} {\bf Proof of #1.}}

\newcommand\setof[1]{\left\{#1\right\}}

\newcommand\propref[1]{Proposition \ref{#1}}
\newcommand\thmref[1]{Theorem \ref{#1}}
\newcommand\lemref[1]{Lemma \ref{#1}}
\newcommand\corref[1]{Corollary \ref{#1}}
\newcommand\defref[1]{Definition \ref{#1}}

\numberwithin{equation}{section}

\input pstricks.tex
\input xy
\xyoption{all}

\begin{document}
\title[
Compact maps and quasi-finite complexes
]%
   {Compact maps and quasi-finite complexes}

\author{M.~Cencelj}
\address{IMFM,
Univerza v Ljubljani,
Jadranska ulica 19,
SI-1111 Ljubljana,
Slovenija }
\email{matija.cencelj@guest.arnes.si}

\author{J.~Dydak}
\address{University of Tennessee, Knoxville, TN 37996, USA}
\email{dydak@math.utk.edu}

\author{J.~Smrekar}
\address{Fakulteta za Matematiko in Fiziko,
Univerza v Ljubljani,
Jadranska ulica 19,
SI-1111 Ljubljana,
Slovenija }
\email{jaka.smrekar@fmf.uni-lj.si}

\author{A.~Vavpeti\v c}
\address{Fakulteta za Matematiko in Fiziko,
Univerza v Ljubljani,
Jadranska ulica 19,
SI-1111 Ljubljana,
Slovenija }

\email{ales.vavpetic@fmf.uni-lj.si}

\author{\v Z.~Virk}
\address{Fakulteta za Matematiko in Fiziko,
Univerza v Ljubljani,
Jadranska ulica 19,
SI-1111 Ljubljana,
Slovenija }

\email{ziga.virk@fmf.uni-lj.si}

\date{ \today
}


\thanks{ 
}

\begin{abstract}
The simplest condition characterizing quasi-finite CW complexes $K$ is the implication
$X\tau_h K\implies \beta(X)\tau K$ for all paracompact spaces $X$. Here are the main results of the paper:
\begin{Thm}  
If $\{K_s\}_{s\in S}$ is a family of pointed quasi-finite complexes, then their wedge
$\bigvee\limits_{s\in S}K_s$ is quasi-finite.
\end{Thm}
\begin{Thm} If $K_1$ and $K_2$ are quasi-finite countable CW complexes, then their join
$K_1\ast K_2$ is quasi-finite.
\end{Thm}

\begin{Thm}
For every quasi-finite CW complex $K$ there is a family $\{K_s\}_{s\in S}$ of countable CW complexes
such that $\bigvee\limits_{s\in S} K_s$ is quasi-finite and is equivalent, over the class of paracompact
spaces, to $K$.
\end{Thm}

\begin{Thm}
Two quasi-finite CW complexes $K$ and $L$ are equivalent over the class of paracompact
spaces if and only if they are equivalent over the class of compact metric spaces.
\end{Thm}

Quasi-finite CW complexes lead naturally to the concept of $X\tau {\mathcal F}$, where
${\mathcal F}$ is a family of maps between CW complexes. We generalize some well-known results
of extension theory using that concept.
\end{abstract}

\maketitle

\medskip
\medskip
\tableofcontents

%
%

\section{Introduction}

\begin{Def}\label{CompactAndLocCompactMapsDef}
A {\it compact map} is a continuous function $f\colon X\to Y$ such that $f(X)$ is contained in a compact
subset of $Y$. If for each $x\in X$ there is a neighborhood $U_x$ such that $f|_{U_x}\colon U_x\to Y$ is
compact then we say that $f$ is {\it locally compact}.
\end{Def}

\par The notation $K\in AE(X)$ or $X\tau K$ means that any continuous map
map $f\colon A\to K$, with $A$ closed in $X$,
extends over $X$. Also, $K\in AE_h(X)$ or $X\tau_h K$ indicates that any map $f\colon A\to K$, with $A$ closed in $X$,
extends over $X$ up to homotopy. Employing the notation of \cite{Dy2} we use $K\in AE_{lc}(X)$ (or $X\tau_{lc} K$)
to mean that any locally compact map defined on a closed subset of $X$ extends over $X$ to a locally compact map.

\begin{Def}\label{AEForCompactMaps}
A CW complex $K$ is an {\it absolute extensor of $X$ with respect to compact maps}
(notation $K\in AE_{cm}(X)$ or $X\tau_{cm} K$) if every compact map $f\colon A\to K$,
$A$ closed in $X$, extends to a compact map $g\colon X\to K$.
\end{Def}

\par It follows from (4) of Proposition \ref{ResultsFromDy2} that for $X$ a paracompact space
and $K$ a CW complex $X\tau_hK$ and $X\tau_{lc}K$ are equivalent. Furthermore, also $X\tau_hK$
and $X\tau K$ are equivalent if either $X$ is, in addition, metrizable or locally compact (apply
(2) of the same Proposition), or $K$ is locally finite (apply (4)).

This paper is devoted to the question of extendability of compact maps and the main problem
we are interested in is the following:

\begin{Problem}\label{CharOfCMExtensors} 
Given a class $\CCC$ of paracompact spaces characterize CW complexes $K$
such that $X\tau_h K$ and $X\tau_{cm} K$ are equivalent for all $X\in\CCC$.
\end{Problem}

\par It turns out that $X\tau_{cm} K$ is equivalent to $\beta(X)\tau K$;
see the more general Proposition \ref{AEForCompactMapsAndCechStoneMap} below.

\begin{Prop}\label{AEForCompactMapsAndCechStone} 
If $X$ is paracompact and $K$ is a CW complex, then $X\tau_{cm}K$
if and only if $\beta(X)\tau K\in AE$. \edokaz
\end{Prop}

For separable metric spaces the following is true (see Section \ref{proofs} for the proof).

\begin{Prop} \label{BetaXTauKImpliesXTauK}
If $X$ is separable metric and $K$ is a CW complex, then $X\tau_{cm}K$
(equivalently, $\beta(X)\tau K$) implies $X\tau K$.
\end{Prop}

\begin{Problem}\label{Does BetaXTauKImplyXTauKQ} 
Let $X$ be paracompact and let $K$ be a CW complex.
Does $X\tau_{cm}K$ (equivalently, $\beta(X)\tau K$) imply $X\tau_h K$?
\end{Problem}

\begin{Def}\label{DefOfWeaklyQF} Let $\CCC$ be a class of paracompact spaces.
A CW complex $K$ such that $X\in \CCC$ and $X\tau_h K$ imply $X\tau_{cm} K$ is called
{\it weakly quasi-finite with respect to $\CCC$} or simply {weakly \it $\CCC$-quasi-finite}.
\end{Def}

In view of \ref{AEForCompactMapsAndCechStone} Problem
\ref{CharOfCMExtensors} is related to the classical question in any dimension theory:
Does \v Cech-Stone compactification preserve dimension? It is so in the theory
of covering dimension (see  \cite{En$_1$} or  \cite{E$_2$}). In cohomological dimension theory it
was a difficult problem (see \cite{Ku}), finally solved in the negative due to a counterexample of
Dranishnikov \cite{DR2} which was followed by those of Dydak and Walsh \cite{DW}, and Karinski \cite{Kar}.
Dydak and Mogilski \cite{DM} were the first to see a connection between existence of extension dimension
preserving compactifications and $K$-invertible maps (see \ref{InvertibleDef} for a more general definition
of invertibility - for the class of all spaces $X$ with $X\tau K$ Definition \ref{InvertibleDef} corresponds
to $K$-invertibility).

Chigogidze \cite{Ch1} established that not only there is a connection but, in fact, the two
notions are equivalent.

\begin{Thm}[Chigogidze] \label{ChigEquivThm}
For a countable simplicial complex $K$ the following conditions are equivalent:
\begin{itemize}
\item[1.] 	If $X$ is any space then $X\tau K$ if and only if $\beta(X)\tau K$.
\item[2.]	There exist a metrizable compactum $E$ with $E\tau K$ and
		a $K$-invertible map $p\colon E\to I^\omega$ onto the Hilbert cube.
\end{itemize}
\end{Thm}

Karasev \cite{K} gave an intrinsic characterization of countable complexes $P$ satisfying \ref{ChigEquivThm} and
called them {\it quasi-finite complexes}. His definition was later generalized by Karasev-Valov \cite{KV1} as follows:

\begin{Def}\label{QuasiFiniteDef}
A CW complex $K$ is called {\it quasi-finite} if there is a function $e$ from the family
of all finite subcomplexes of $K$ to itself satisfying the following properties:

\begin{itemize}
\item[a.]	$M\subset e(M)$ for all $M$.
\item[b.]	If $X$ is a Tychonoff space, $K$ is an absolute extensor for $X$
		in the sense of \cite{Ch2}, and $A$ is any closed subset of $X$
		then every map $f\colon A\to M$ extends over $X$ to a map $g\colon X\to e(M)$.
\end{itemize}
\end{Def}

Obviously our definition \ref{DefOfWeaklyQF} of quasi-finiteness is weaker than \ref{QuasiFiniteDef}.
However, later on we show that the two are equivalent.


In subsequent sections we discuss the relationship between properties $X\tau_h K$ and $\beta(X)\tau K$,
under different hypotheses on $X$ and $K$. The main focus is on quasi-finite CW complexes $K$.

In Section \ref{families} we generalize the concept $X\tau_h K$ to the concept $X\tau_h\FFF$ where
$\FFF$ is an arbitrary family of maps between CW complexes. It turns out that this concept is natural
for extension theory. We introduce various notions of quasi-finite families and discuss relations among
them. We also generalize the notion of a $K$-invertible map (see \cite{KV1}) for a CW complex $K$ to
an $\FFF$-invertible map where $\FFF$ is a family of maps between CW complexes. We discuss relations
between existence of invertible maps and quasi-finite families.

In Section \ref{qfcomplexes} we apply our results to quasi-finite CW complexes, ie the case when the
family $\FFF$ solely consists of the identity map on a CW complex, and deduce several geometric properties
of quasi-finite CW complexes.

In Section \ref{qfdim} we prove the existence of ``quasi-finite dimension'' for a given paracompact space $X$.

Section \ref{proofs} contains most proofs, and Section \ref{op} gives a short list of problems -- open, to the
best of our knowledge.

We would like to thank Leonard Rubin for many comments and suggestions which were incorporated into our paper.


\section{Extension properties of families of maps}\label{families}

In \cite{DrR1} Dranishnikov and Repov\v s introduced the concept of $X\tau(L,K)$,
where $K\subset L$. It means that any map $f\colon A\to K$, with $A$ closed in $X$,
extends to $g\colon X\to L$. Rubin (see \cite{Rubin}, Definition 6.1) says in that
case that the pair $(L,K)$ is $X$-connected. In \cite{DR3} Dranishnikov considers the
concept of $X\tau i$, where $i\colon K\to L$ is a map. This means that for any map
$f\colon A\to K$ defined on the closed subset $A$ of $X$ the composite $if\colon A\to L$
extends over $X$. We depict this in the following diagram.
\begin{equation*}\tag{$*$}
\begin{CD} A @>f>> K \\ @V{\text{incl}}VV @VV{i}V \\ X @>F>> L
\end{CD}
\end{equation*}

We can also define $X\tau_hi$ by demanding the existence of $F$ so that the diagram ($*$)
commutes up to homotopy. Furthermore, we consider the following notions:
\par
We say that $X\tau_{cm}i$ (respectively $X\tau_{lc}i$) if for any compact (respectively
locally compact) map $f$ there exists a compact (respectively locally compact) map $F$
making the diagram ($*$) commutative. The set $A$ is required to be closed in $X$.

We generalize the above notions as follows.

\begin{Def}\label{NewExtensionDef}
Suppose ${\mathcal F}$ is a family of maps between CW complexes.
Given a paracompact space $X$ we say that $X\tau {\mathcal F}$ or ${\mathcal F}\in AE(X)$
(respectively $X\tau_h\FFF$, $X\tau_{cm}\FFF$, $X\tau_{lc}\FFF$) if $X\tau i$ (respectively
$X\tau_h i$, $X\tau_{cm}i$, $X\tau_{lc}i$) for all $i\in\FFF$.
\end{Def}

As in the case of absolute extensors, the motivating problem is the relation between the properties
$X\tau_h\FFF$ and $\beta(X)\tau\FFF$. Here the \v Cech-Stone compactification of $X$ can be replaced
by $X$ if we restrict to compact maps:

\begin{Prop}\label{AEForCompactMapsAndCechStoneMap} 
Let $X$ be paracompact and let $\FFF$ be a family of maps between CW complexes.
Then $X\tau_{cm}\FFF$ if and only if $\beta(X)\tau \FFF$.
\end{Prop}

If $\FFF$ is a family of maps between finite CW complexes, the following is true.

\begin{Prop} \label{quasi-finite}
Let $X$ be a paracompact space and $\FFF$ a family of maps between
finite (or finitely dominated) CW complexes. Then $X\tau_h\FFF$ if and only if $\beta(X)\tau\FFF$.
\end{Prop}

\begin{Question}\label{TauImpliesBetaTau}
Find necessary and sufficient conditions for $X\tau_h\FFF$ to imply $\beta(X)\tau\FFF$ for all $X\in\CCC$.
\end{Question}

We need the notion of {\it extensionally equivalent} families.

%
%

\subsection{Extensional equivalence of families}

\begin{Def}\label{EquivOfFamiliesDef}
Suppose ${\mathcal C}$ is a class of paracompact spaces. Two families ${\mathcal F}$ and ${\mathcal G}$
of maps between CW complexes are {\it equivalent over ${\mathcal C}$} (or simply ${\mathcal C}$-equivalent)
if $X\tau_h {\mathcal F}$ is equivalent to $X\tau_h {\mathcal G}$ for all $X\in {\mathcal C}$.
\par More generally, if $X\tau_h\FFF$ implies $X\tau_h\GGG$ whenever $X\in\CCC$, we say that
{\it $\FFF$ is extensionally smaller than $\GGG$ over $\CCC$} and denote it $\FFF\leq\GGG$ over
$\CCC$ or simply $\FFF\leq_\CCC\GGG$.
\end{Def}

Extensionally, we can always replace a family of maps by a single map.

\begin{Thm}\label{family_equivalent_to_map}
Let $\FFF=\{i_s\colon K_s\to L_s\,\vert\,s\in S\}$ be a family of maps between CW complexes,
and form the obvious map $i\colon\vee_sK_s\to\vee_sL_s$. Then $\FFF$ is equivalent to $\{i\}$
over the class of paracompact spaces.
\end{Thm}

\begin{Rem}
Although it follows from Theorem \ref{family_equivalent_to_map} that every family of maps
between CW complexes is equivalent over the class of paracompact spaces to a single map between
CW complexes, it makes sense to consider families of maps. Namely, in light of Proposition \ref{quasi-finite}
we are interested in precisely those maps (or families of maps) which are equivalent to a family
of maps between {\bf finite} CW complexes.
\end{Rem}

The following result shows that extension properties of a family with respect to
the class of compact Hausdorff spaces are the same as those with respect to compact metric spaces.

\begin{Thm} \label{EquivOverCompactaImpliesOverCompact}
Suppose ${\mathcal F}$ and ${\mathcal G}$ are families of maps between
CW complexes. If ${\mathcal F}\leq{\mathcal G}$ over
compact metric spaces, then ${\mathcal F}\leq{\mathcal G}$ over compact Hausdorff spaces.
\end{Thm}

The class of compact metric spaces is important enough to justify the following definition.

\begin{Def}
Let $\CCC$ be a class of paracompact spaces. We say that $\CCC$ is {\it rich} if it contains
all compact metric spaces.
\end{Def}

%
%

\subsection{The notions of quasi-finiteness}

For the purpose of answering Question \ref{TauImpliesBetaTau} we introduce the following notions.

\begin{Def}\label{QFNotions}\label{EquFiniteOfFamiliesDef}\label{StronglyquFiniteOfFamiliesDef}
Suppose $\CCC$ is a class of paracompact spaces and $\FFF$ is a family of maps between CW complexes.
\begin{enumerate}
\item	The family $\FFF$ is {\it weakly quasi-finite over ${\mathcal C}$}
	(or simply weakly ${\mathcal C}$-quasi-finite) if
	$X\tau_h {\mathcal F}$ implies $\beta(X)\tau {\mathcal F}$ whenever $X\in\CCC$
\item	The family $\FFF$ is {\it quasi-finite over ${\mathcal C}$} (or simply ${\mathcal C}$-quasi-finite) if
	for each member $i\colon K\to L$ of $\FFF$ and for each finite subcomplex $M$ of $K$
	there is a finite subcomplex $e_i(M)$ of $L$ containing the image $i(M)$ such that
	all restrictions $j_{i,M}\colon M\to e_i(M)$ of $i$ satisfy $X\tau j_{i,M}$
	if $X\in\CCC$ satisfies $X\tau_h\FFF$.
\item	The family $\FFF$ is {\it strongly quasi-finite over ${\mathcal C}$}
	(or strongly ${\mathcal C}$-quasi-finite) if ${\mathcal F}$ is
	${\mathcal C}$-equivalent to a family ${\mathcal G}$ of maps between finite CW complexes.
\end{enumerate}
\end{Def}

\begin{Def}\label{QFCWcx}
Suppose that $\CCC$ is a class of paracompact spaces and let $K$ be a CW complex.
Then we say that $K$ is {\it weakly $\CCC$-quasi-finite} or {\it $\CCC$-quasi-finite} or
{\it strongly $\CCC$-quasi-finite} if the property in question holds for the singleton family
$\setof{\id\colon K\to K}$.
\end{Def}

It is obvious that the notion of a strongly quasi-finite family $\FFF$ with respect to an arbitrary class is
actually well defined for the equivalence class of $\FFF$. For the notion of a weakly quasi-finite family
the same is true if the reference class of paracompact spaces is rich:

\begin{Prop}\label{EquivToQFIsQF} 
Let $\CCC$ be a rich class of paracompact spaces, and suppose that families
$\FFF$ and $\GGG$ of maps between CW complexes are equivalent over ${\mathcal C}$.
If $\FFF$ is weakly ${\mathcal C}$-quasi-finite, then so is $\GGG$.
\end{Prop}
\dokaz 
By \thmref{EquivOverCompactaImpliesOverCompact}, the families $\FFF$ and $\GGG$ are equivalent
over the class of all compact Hausdorff spaces. 
Suppose $X\in {\mathcal C}$ and $X\tau_h \GGG$. Hence $X\tau_h \FFF$ 
which implies $\beta(X)\tau \FFF$. Hence $\beta(X)\tau \GGG$
and $L$ is weakly quasi-finite with respect to ${\mathcal C}$.
\edokaz

\begin{Def}
Let $\CCC$ be a class of paracompact spaces, let $\FFF$ be a family of maps between CW complexes,
and let $\alpha$ be an infinite cardinal number.
\begin{itemize}
\item	The class $\CCC$ is {\it $\alpha$-saturated} if it is closed under disjoint unions of at most $\alpha$ members.
	If $\CCC$ is $\alpha$-saturated with respect to all $\alpha$ then $\CCC$ is called {\it saturated}.
\item	The family $\FFF$ is {\it range-dominated by $\alpha$} is the ranges of all maps in $\FFF$ have at most $\alpha$
	cells (or are homotopy dominated by such complexes).
\end{itemize}
Let $\CCC$ be any class of paracompact spaces.
Let $\Bar\CCC^\alpha$ denote all possible disjoint unions of at most $\alpha$ members of $\CCC$,
and let $\Bar\CCC$ denote all possible disjoint unions of members of $\CCC$.
Then $\Bar\CCC^\alpha$ is evidently $\alpha$-saturated while $\Bar\CCC$ is saturated.

We call $\Bar\CCC^\alpha$ and $\Bar\CCC$ the {\it $\alpha$-saturation} and the {\it saturation} of $\CCC$,
respectively.
\end{Def}

We note the evident

\begin{Prop}\label{saturation}
Let $\FFF$ be weakly quasi-finite (respectively strongly quasi-finite) with respect to $\CCC$.
Then $\FFF$ is actually weakly quasi-finite (respectively strongly quasi-finite) with respect to $\Bar\CCC$.\edokaz
\end{Prop}

\begin{Thm}\label{TheQFThm}\label{QFImpliesWQF}\label{ExistenceOfSupComplexForFamilies}\label{StronglyQFAreWQF}
Let $\CCC$ be a class of paracompact spaces and let $\FFF$ be a family of maps between CW complexes.
\begin{enumerate}
\item	If $\FFF$ is $\CCC$-quasi-finite then $\FFF$ is weakly $\CCC$-quasi-finite.
\item	If $\FFF$ is range-dominated by the infinite cardinal number $\alpha$ and $\CCC$ is $\alpha$-saturated
	then $\FFF$ is $\CCC$-quasi-finite if and only if $\FFF$ is weakly $\CCC$-quasi-finite.
\item	If $\FFF$ is strongly $\CCC$-quasi-finite and $\CCC$ is rich then $\FFF$ is weakly $\CCC$-quasi-finite.
\end{enumerate}
\end{Thm}




\begin{Thm}\label{EquivOverCompAreEquiv}
Let $\CCC$ be a rich class of paracompact spaces. Two strongly $\CCC$-quasi-finite
families are equivalent over $\CCC$ if and only if they are equivalent over the class
of compact metric spaces.
\end{Thm}

The following theorem partially answers our motivating problems.

\begin{Thm} \label{CharOfQFFamilies}
Let $\CCC$ be a class of paracompact spaces and suppose that the family $\FFF$ is $\CCC$-quasi-finite.
\begin{enumerate}
\item	If $\FFF$ is strongly $\CCC$-quasi-finite and $\CCC$ is rich then 
	$X\tau_h\FFF$ is equivalent to $\beta(X)\tau\FFF$ whenever $X\in\CCC$.
\item	If $X\tau_h\FFF$ is equivalent to $\beta(X)\tau\FFF$ whenever $X\in\CCC$
	then $\FFF$ is strongly $\CCC$-quasi-finite.
\end{enumerate}
\end{Thm}

%
%

\subsection{Invertible maps}

\begin{Def}\label{InvertibleDef}
Let $\CCC$ be a class of paracompact spaces and let $\FFF$ be a family of maps between CW complexes.
The map $p\colon E\to B$ is {\it $\FFF$-invertible} with respect to $\CCC$ if any map $f\colon X\to B$
with $X\tau_h\FFF$ and $X\in\CCC$ has a lift $g\colon X\to E$.
\end{Def}

\begin{Rem}\label{rem_inv_sat}
Note that for any class $\CCC$, a given map $p\colon E\to B$ is $\FFF$-invertible with respect to
$\CCC$ if and only if $p$ is $\FFF$-invertible with respect to the saturation $\Bar\CCC$.
\end{Rem}

\begin{Prop}\label{ExistenceOfInvertibleImpliesWQF}
Let $\CCC$ be a class of paracompact spaces. If there is an $\FFF$-invertible map
$p\colon E\to I^\omega$ with respect to $\CCC$ where $E$ is a compact Hausdorff space
and $E\tau \FFF$, then $\FFF$ is weakly $\CCC$-quasi-finite.
\end{Prop} 

Recall that the {\it weight} of a topological space is defined in the following
way. Let $\alpha$ be a cardinal number, and let $Y$ be a topological space.
We say that the weight of $Y$ is less than $\alpha$ (or is at most $\alpha$),
in symbols $w(Y)\leq\alpha$, if the topology on $Y$ admits a basis whose
cardinality is at most $\alpha$.

\begin{Thm}\label{ExistenceOfInvertibleMapsForQFThm}
Let ${\mathcal C}$ be a class of paracompact spaces and let $\FFF$ be a family of maps between CW complexes.
If $\FFF$ is $\CCC$-quasi-finite, then for each infinite cardinal number $\alpha$ there are a compact Hausdorff space $E$
of weight at most $\alpha$ and an $\FFF$-invertible map $p\colon E\to I^\alpha$ with respect to $\CCC$ such that
$E\tau\FFF$.
\end{Thm} 

\begin{Cor}[Karasev-Valov \cite{KV1}, Theorem 2.1]\label{InvExistCor}
If $K$ is a quasi-finite CW complex (with respect to paracompact spaces),
then for any $\alpha$ there are a compact Hausdorff space $E$ of weight
at most $\alpha$ and a $K$-invertible map $p\colon E\to I^\alpha$ such that $E\tau K$. \edokaz
\end{Cor}

%
%

\subsection{Geometric realization of families of maps}

\begin{Def}\label{GeomRealizationDef}
Let $\FFF$ be a family of maps between CW complexes and let $\CCC$ be a class of paracompact spaces.
A CW complex $K$ is {\it a geometric $\CCC$-realization} of $\FFF$ if $X\tau_h\FFF$ is equivalent to $X\tau_h K$
whenever $X\in\CCC$.
\end{Def}

Note that a CW complex $K$ is strongly $\CCC$-quasi-finite if and only if $K$ is a geometric $\CCC$-realization
of a family of maps between finite CW complexes.

\begin{Prop}\label{RealizationOfUnionOfFamilies} 
Let $\setof{\FFF_s\,\vert\,s\in S}$ be families of maps between CW complexes indexed by the set $S$.
If $K_s$ is a geometric $\CCC$-realization of $\FFF_s$ for each $s$, then $\bigvee\limits_{s\in S}K_s$
is a geometric $\CCC$-realization of $\bigcup\limits_{s\in S}\FFF_s$.
\end{Prop}

\begin{Problem}\label{RealizationOfJoinOfFamilies} 
Suppose $K_s$ is a geometric realization of ${\mathcal F}_s$ for each $s=1,2$.
Is $K_1\ast K_2$ a geometric realization of ${\mathcal F}_1\ast {\mathcal F}_2$?

By $\FFF_1\ast\FFF_2$ we mean the collection of maps $i_1\ast i_2$ for all $i_1$ and $i_2$
that belong to $\FFF_1$ and $\FFF_2$, respectively. See the proof of \thmref{JoinOfQF} on page
\pageref{ProofJoinOfQF} for an explanation of this definition of $\FFF_1\ast\FFF_2$.

To avoid technical difficulties surrounding the choice of topology on $K_1\ast K_2$
(in general, namely, we have to take the compactly generated refinement of the topological join
to get the natural CW structure) one can start with the case of $K_1$ and $K_2$ countable. In
that case the topological join $K_1\ast K_2$ is already compactly generated.
\end{Problem}

\begin{Problem}\label{ExistenceOfFWithNoQFRealization} 
Is there a family ${\mathcal F}=\{i_s\colon K_s\to L_s\}_{s\in S}$
such that all $K_s$, $L_s$ are finite and ${\mathcal F}$ has no geometric realization?
\end{Problem}

\begin{Example}
Let $p$ be a natural number and let $f\colon K\to K$ be a self-map. If the iterate $f^{p+1}$ is homotopic
to $f$, then $f$ is called {\it a homotopy $p$-idempotent} (see \cite{Smr}). The case of ordinary homotopy
idempotents occurs for $p=1$.

Let $\Tel_f$ denote the infinite mapping telescope of $f$, that is, the quotient space of
$K\times\NN\times[0,1]$ modulo the relations $(x,n,1)\simeq(f(x),n+1,0)$. If $h$ is
a homotopy between the iterates $f$ and $f^{p+1}$, then a map $u\colon\Tel_f\to K$ may be defined
by \[ [x,n,t]\mapsto f^{n(p-1)}h(x,t). \]
Here $[x,n,t]$ denotes the equivalence class of $(x,n,t)$ in $\Tel_f$. There
is an obvious map $d\colon K\to\Tel_f$ defined by $d(x)=[x,0,0]$.

Observe that the composite $ud$ equals $f$. The $p$-idempotent $f$ {\it splits} (see \cite{Smr})
if the composite $du$ is a homotopy equivalence. In this case $(du)^p$ is homotopic to the identity.

By the cellular approximation theorem we can assume $f$ to be cellular. In that case the infinite
mapping telescope $\Tel_f$ admits an obvious CW decomposition.
The following proposition says that $\Tel_f$ is a geometric realization
of $\setof{f}$ if $f$ is a split $p$-idempotent.
\end{Example}

\begin{Prop}\label{HomotopyIdempotentOfFinite} 
Let $f\colon K\to K$ be a cellular homotopy $p$-idempotent and
let $\Tel_f$ be the infinite mapping telescope of $f$.
Let $X$ be any space.
\begin{enumerate}
\item	$X\tau_h\Tel_f$ implies $X\tau_h f$.
\item	If $f$ splits, then $X\tau_h f$ implies $X\tau_h\Tel_f$.
\end{enumerate}
\end{Prop}
\dokaz 
An immediate consequence of Lemma \ref{ReductionIfTargetIsHomotopyDominated}.
\edokaz

\begin{Cor}\label{CorHomotopyIdempotentOfFinite}
Let $f\colon K\to K$ be a cellular homotopy $p$-idempotent on the finite-dimensional CW complex $K$,
that is, a cellular map for which $f^{p+1}$ -- the $(p+1)$st iterate of $f$ -- is homotopic to $f$.
Then the infinite mapping telescope of $f$ is a geometric realization of $\setof{f\colon K\to K}$.
\end{Cor}
\dokaz
Follows from \propref{HomotopyIdempotentOfFinite} together with Corollary 2.4 of \cite{Smr}.
\edokaz

\begin{Example}\label{PuncturedTorusProp} 
Let $T$ be a torus with an open disk removed and let $\partial T$ be the boundary of $T$.
If $K$ is a CW complex such that $K$ geometrically realizes the inclusion $i_T\colon\partial T\to T$,
then $K$ is acyclic.
\end{Example}
\dokaz As in \cite{DrR1} (Theorem 4.8) the $n$-sphere $S^n$, $n\ge 2$, can be
split as $X_1\cup X_2$ where $X_2$ is $0$-dimensional
and $X_1$ is a countable union of compact metric spaces $A_k$ satisfying
$A_k\tau i_T$ for each $k$. Therefore $X_2\tau K$
and, by the union theorem for extension theory (see \cite{Dy3}), it follows that
$S^n\tau S^0\ast K$. This is only possible if $S^0\ast K=\Sigma(K)$ is contractible.
\edokaz

%
%

\section{Geometric properties of quasi-finite CW complexes}\label{qfcomplexes}

Here is the main property of quasi-finite CW complexes.

\begin{Thm} \label{CharOfQFViaSepMetric}
Let $K$ be an infinite CW complex of cardinality $\alpha$ and let $\CCC$ be a class
of paracompact spaces. Consider the following statements.
\begin{enumerate}
\item $K$ is strongly $\CCC$-quasi-finite.
\item $K$ is $\CCC$-quasi-finite.
\item $K$ is weakly $\CCC$-quasi-finite.
\end{enumerate}
\par The implication (2)$\implies$(1) always holds while (1)$\implies$(3) holds if $\CCC$ is rich.
\par The implication (3)$\implies$(2) holds if $\CCC$ is $\alpha$-saturated.
\par In particular, (1),(2),(3) are equivalent if $\CCC$ is a rich $\alpha$-saturated class of paracompact spaces.
\end{Thm}

\begin{Rem}
For theorem 7.6 of \cite{Rubin}, the `cycle' of implications (d)$\implies$(e)$\implies$(f)$\implies$(c)$\implies$(d)
can also be deduced from our results in the following way. (d)$\implies$(f) follows from
\thmref{ExistenceOfInvertibleMapsForQFThm}, (f)$\implies$(e) is a tautology, (e)$\implies$(c) follows
from Remark \ref{rem_inv_sat} and \propref{ExistenceOfInvertibleImpliesWQF}, and (c)$\implies$(d) follows from
\thmref{CharOfQFViaSepMetric}.

The equivalence (b)$\iff$(c) of \cite{Rubin} shows, in our language, that a CW complex $K$
is weakly quasi-finite with respect to the $\aleph_0$-saturation of compact metric spaces
if and only if $K$ is weakly quasi-finite with respect to the (total) saturation of the class
of compact metric spaces. We consider that result very interesting.
\end{Rem}

In light of \propref{saturation} we immediately deduce

\begin{Cor} \label{CorCharOfQFViaSepMetric}
Let $K$ be an infinite CW complex and $\CCC$ a rich class of paracompact spaces.
Then $K$ is $\CCC$-quasi-finite if and only if it is strongly $\CCC$-quasi-finite.\edokaz
\end{Cor}

If $K$ is countable then the notion of quasi-finiteness in our \defref{DefOfWeaklyQF}
is equivalent that of \cite{K}.

\begin{Cor}\label{ModOfOriginalDefOfQF}
If $K$ is a countable CW complex, then the following are equivalent:
\begin{itemize}
\item[a.] $K$ is quasi-finite with respect to Polish spaces.
\item[b.] $K$ is quasi-finite with respect to paracompact spaces.
\end{itemize}
\end{Cor}

\begin{Cor}\label{QFAndCompactificationsOfSeparable}
If $K$ is a countable CW complex, then the following are equivalent:
\begin{itemize}
\item[a.] $K$ is quasi-finite with respect to Polish spaces.
\item[b.] For each metric compactification $\nu(X)$ of a separable metric space $X$
satisfying $X\tau K$ there is a metric compactification $\gamma(X)$ of $X$
such that $\gamma(X)\tau K$ and $\gamma(X)\ge \nu(X)$.
\end{itemize}
\end{Cor}

\begin{Thm}\label{WedgeOfQFIsQF}
If $\{K_s\}_{s\in S}$ is a family of pointed ${\mathcal C}$-quasi-finite 
(respectively, strongly or weakly $\CCC$-quasi-finite) complexes, then their wedge
$K=\bigvee\limits_{s\in S}K_s$ is ${\mathcal C}$-quasi-finite (respectively, strongly or weakly $\CCC$-quasi-finite).
\end{Thm}

\begin{Cor}\label{ExampleOfWedgeNonQF}
The wedge $\bigvee\limits_{p\text{ prime}}K(Z_{(p)},1)$ is not quasi-finite.
\end{Cor}
\dokaz
By the well-known First Theorem of Bockstein the CW complex $\bigvee\limits_{p\text{ prime}}K(Z_{(p)},1)$
is equivalent, over the class of compact metric spaces, to $S^1$. However, Dranishnikov-Repov\v s-Shchepin
\cite{DrRS} proved the existence of a separable metric space $X$ of dimension $2$ such that
$\dim_{Z_{(p)}}(X)=1$ for all prime $p$. Thus, $\bigvee\limits_{p\text{ prime}}K(Z_{(p)},1)$ is not
equivalent to $S^1$ over the class of separable metric spaces and is therefore not quasi-finite.
\edokaz

In \cite{DR1} Dranishnikov proved that every CW complex is equivalent, over the class of compact Hausdorff spaces,
to a wedge of countable CW complexes. Our next result is a variation of that theorem.

\begin{Thm}\label{EveryQFIsWedgeOfCount}
For every $\CCC$-quasi-finite CW complex $K$ there is a family $\{K_s\}_{s\in S}$ of countable CW complexes
so that $\bigvee\limits_{s\in S}K_s$ is $\CCC$-quasi-finite and $\CCC$-equivalent to $K$.
\end{Thm}

\begin{Thm}\label{UnionThmForMaps}
Let $j_i\colon K_i\to L_i$, $i=1,2$, be maps between CW complexes. Suppose $X$ is a separable metric space
and $X_1$, $X_2$ are subsets of $X$. If $Y_1\tau j_1$ for every subset $Y_1$ of $X_1$
and $Y_2\tau j_2$ for every subset $Y_2$ of $X_2$, then $X_1\cup X_2\tau j_1\ast j_2$.
\end{Thm}

\begin{Thm}\label{JoinOfQF}  If $K_1$ and $K_2$ are quasi-finite countable complexes, then their join
$K_1\ast K_2$ is quasi-finite.
\end{Thm}

%
%

\section{Quasi-finite dimension}\label{qfdim}

For each paracompact space $X$ consider all quasi-finite complexes $K$ such that $X\tau_h K$.
That class has an initial element $K_X$ (with respect to the relation $K\leq L$ over paracompact spaces)
described as follows: consider all countable CW complexes $K_s$, $s\in S_X$, that
appear in a decomposition guaranteed by \thmref{EveryQFIsWedgeOfCount}
of a quasi-finite CW complex $K$ with $X\tau_h K$. To make sure that $S_X$
is a set, we choose only one representative $K_s$ within its homotopy type.
(There is only a set of distinct homotopy types of countable CW complexes.)
The wedge $K_X$ of those $K_s$ is quasi-finite. Indeed, $K_X$ is equivalent
to the wedge of $S_X$ copies of itself and that wedge is equivalent to a wedge of quasi-finite complexes.
By \thmref{WedgeOfQFIsQF} $K_X$ is quasi-finite. Also, $K_X\leq K$ for all quasi-finite $K$
satisfying $X\tau_h K$ which means that $K_X$ is indeed an initial element.
It is called the {\it quasi-finite dimension} of $X$ and is denoted by $\dim_{QF}(X)$.

\begin{Thm} \label{QFDimRealizedByCompacta}
For every quasi-finite CW complex $K$ there is a compactum $X$ such that
$\dim_{QF}(X)=K$ and the extension dimension of $X$ equals $K$.
\end{Thm}

\begin{Prop}\label{PropsOfQFDim}
\begin{enumerate}
\item	If $X$ is any paracompact space then $\dim_{QF}(X)=\dim_{QF}(\beta(X))$.
\item	If $X$ is a separable metric space then there exists a metric compactification
	$c(X)$ of $X$ such that $\dim_{QF}(X)=\dim_{QF}(c(X))$.
\item	If $Y$ is a separable metric space and $X\subset Y$, then there
	is a $G_\delta$-subset $X'$ of $Y$ containing $X$ such that $\dim_{QF}(X)=\dim_{QF}(X')$.
\end{enumerate}
\end{Prop}

An alternative way to define quasi-finite dimension is
to fix a countable family $\{M_i\}_{i=1}^\infty$ of finite CW complexes
such that any finite CW complex is homotopy equivalent to
one of them and assign to each paracompact space $X$ the family
of homotopy classes of maps $f\colon M_i\to M_j$ such that there is a quasi-finite CW complex $K_f$ satisfying
$\{id\colon K_f\to K_f\}\leq f$ over paracompact spaces and $X\tau_h K_f$.
The advantage of this approach is that $\dim_{QF}(X)$
would be a countable object. Note that the family of such $f$ would be equivalent
to $\{id\colon K\to K\}$, where $K$ is the wedge of all $K_f$ and $K$ is quasi-finite.
In that sense the two approaches are equivalent.

%
%

\section{Proofs}\label{proofs}


\dokazp{\propref{BetaXTauKImpliesXTauK}}
Consider $X$ as a subset of the Hilbert cube $I^\omega$ and let $f\colon \beta(X)\to I^\omega$
be an extension of the inclusion $X\to I^\omega$. By Theorem 1.1 of \cite{LRS} there is
a factorization $f=p\circ g$ through $g\colon\beta(X)\to Y$ where $Y$ is compact metric and $Y\tau K$.
Note that $Y$ is a metric compactification of $X$. Corollary 3.7 of \cite{ivansic-rubin} implies $X\tau K$.
\edokaz

For convenience we record a few results concerning the relations between
$X\tau\FFF$, $X\tau_h\FFF$ and $X\tau_{lc}\FFF$ implied by \cite{Dy2}.

\begin{Prop}\label{ResultsFromDy2}
\begin{enumerate}
\item	If $f$ is a locally compact map then $\psi f\varphi$ is locally compact
	for any maps $\varphi,\psi$ for which the composition makes sense.
\item	If $X$ is first countable or locally compact Hausdorff then any map
	from $X$ to a CW complex is locally compact. In particular for such $X$
	and $i$ any map between CW complexes, $X\tau i$ if and only if $X\tau_hi$.
\item	If $X$ is paracompact and $f\colon A\to K$ is a locally compact map into the
	CW complex $K$, defined on the closed subset $A$ of $X$, then $f$ extends over $X$
	up to homotopy if and only if it extends over $X$ to a locally compact map.
\item	Let $X$ be paracompact and let $i\colon K\to L$ be a map between CW complexes.
	\begin{itemize}
\item[a.] $X\tau_{h}i$ if and only if $X\tau_{lc}i$.
\item[b.]  If, in addition, $i$ is locally compact,
	then $X\tau_{h}i$ if and only if $X\tau i$ as well.

\end{itemize}
\end{enumerate}
\end{Prop}
\dokaz
Statement (1) is trivial, while statement (2) follows from Corollary 5.4 of \cite{Dy2}.
Statement (3) follows from Corollary 2.13 of \cite{Dy2}, while (4) is an immediate
consequence of (3).
\edokaz

\dokazp{Proposition \ref{AEForCompactMapsAndCechStoneMap}}
It suffices to consider the case of $\FFF$ consisting of one map $i\colon K\to L$.
Suppose $X\tau_{cm}i$ and let $f\colon A\to K$ be a map defined on the closed subset $A$ of $\beta(X)$.
The image of $f$ is contained in a finite subcomplex $K'$ of $K$. By Corollary 2.11 of \cite{Dy2}, the
map $f\colon A\to K'$ extends (strictly) over a neighborhood $N$ of $A$ in $\beta(X)$. We may assume
that $N$ is closed in $\beta(X)$ and abuse notation to let $f\colon N\to K'$ denote the extension.
Since $f|_{N\cap X}\colon N\cap X\to K'$ is a compact map, the composite $i\circ f\vert_{N\cap X}$ admits
a compact extension $h\colon X\to L$, that is, the image of $h$ is contained in a finite subcomplex
$L'$ of $L$. Viewing $L'$ as a subset of the Hilbert cube we may extend $h$ to a map $F\colon\beta(X)\to L'$.
Since $X\cap\Int(N)$ is dense in $\Int(N)$ it follows that $F|_{\Int(N)}=i\circ f|_{\Int(N)}$.
In particular, $F|_A=i\circ f$.
\par Suppose now that $\beta(X)\tau i$ and let $f\colon A\to K$ be a compact map defined on the closed
subset $A$ of $X$. Certainly $f$ extends over $\beta(A)$ which is closed in $\beta(X)$. Consequently
$i\circ f$ extends to a map $\beta(X)\to L$, and the restriction of that extension to $X$ is a compact
extension of $i\circ f$.
\edokaz

The following is a generalization of the Marde\v si\'c
(or Levin, Rubin, Schapiro, see \cite{LRS}) factorization theorem.

\begin{Thm} \label{FactorizationThmForFamily}
Let $\alpha$ be an infinite cardinal.
If $f\colon X\to Y$ is a map of compact Hausdorff spaces with $w(Y)\leq\alpha$,
then $f$ factors as $f=p\circ g$, where $g\colon X\to Z$ is surjective,
$Z$ is compact, $w(Z)\leq\alpha$, and for any family $\FFF$ of maps between CW complexes
$X\tau\FFF$ implies $Z\tau\FFF$.
\end{Thm}
\dokaz Theorem 1.1 of \cite{LRS} implies that $f$ factors as $f=p\circ g$,
where $g\colon X\to Z$ is surjective, $Z$ is compact metric and for any map $f_0\colon A\to M$,
with $A$ closed in $Z$ and $M$ a CW complex, $f_0$ extends over $Z$ if
$f_0\circ g\colon g^{-1}(A)\to M$ extends over $X$.

Let $i\colon K\to L$ be a member of $\FFF$.
Suppose $h\colon A\to K$ is a map from a closed subset $A$ of $Z$.
Put $f_0=i\circ h\colon A\to L$. If $X\tau i$, an extension
$F\colon g^{-1}(A)\to L$ of $f_0\circ g\vert_{g^{-1}(A)}$ does exist.
Therefore $f_0$ extends over $Z$. This shows $Z\tau i$.
\edokaz

%
%

\begin{Lem}\label{ReductionIfTargetIsHomotopyDominated}
Let $f\colon K\to L$ be a map with $L$ and $K$ homotopy dominated by
$N$ and $M$ respectively (i.e., there are maps $d\colon N\to L$, $u\colon L\to N$, $d'\colon M\to K$, $u'\colon K\to M$
with $du$ homotopic to $\id_L$ and $d'u'$ homotopic to $\id_K$).
For any space $X$ the following are equivalent:
\begin{itemize}
\item[i.] $X\tau_h f$,
\item[ii.] $X\tau_h(uf)$,
\item[iii.] $X\tau_h(fd')$.
\end{itemize}
\end{Lem}

\dokaz
Note that if $\lambda\colon L\to L'$ is any map, $X\tau_h f$ implies $X\tau_h(\lambda f)$. Dually, if $\kappa\colon K'\to K$
is any map, $X\tau_h f$ implies $X\tau_h(f\kappa)$. This establishes implications ({\bf i})$\implies$({\bf ii}) and
({\bf i})$\implies$({\bf iii}).

The same argument settles the reverse implications as $f$ is homotopic to $d(uf)$ and $(fd')u'$.
\hfill $\blacksquare$

\begin{Lem} \label{CechStoneCompLem}
Let $X$ be a paracompact space and let $i\colon K\to L$ be a map between CW complexes.
If $L$ is homotopy dominated by a finite CW complex, then $X\tau_h i\implies\beta(X)\tau i$.

Consequently if ${\mathcal F}$ is a family of maps with ranges finitely dominated,
$X\tau_h {\mathcal F}$ implies $\beta(X)\tau  {\mathcal F}$.
\end{Lem}
\dokaz
First note that since $\beta X$ is compact Hausdorff, $\beta X\tau i$ is equivalent to $\beta X\tau_h i$,
by (2) of Proposition \ref{ResultsFromDy2}. Therefore we may assume that $L$ is in fact finite, by ({\bf ii})
of Lemma \ref{ReductionIfTargetIsHomotopyDominated}.

By (4) of Proposition \ref{ResultsFromDy2}, $X\tau_{lc}i$ hence $X\tau_{cm}i$ as $L$ is compact. Thus
$\beta(X)\tau i$ by Proposition \ref{AEForCompactMapsAndCechStoneMap}.
\hfill $\blacksquare$

\begin{Lem} Let $A$ be a closed subset in $\beta X$. Then $\beta X=\beta(X \cup A)$.
\end{Lem}
\dokaz
Take $g \colon X \cup A \to I$. We claim that there exists a unique extension to $\beta X \to I$.
Certainly there exists an extension $h\colon \beta X\to I$ of $g\vert_X\colon X\to I$. As $X$ is
dense in $X\cup A$ and $I$ is Hausdorff, $h\vert_{X\cup A}=g$.
\hfill $\blacksquare$

\begin{Lem} \label{ReverseOfCechStoneCompLem}
Let $X$ be a paracompact space and let $i\colon K\to L$ be a map between CW complexes.
If $K$ is homotopy dominated by a finite CW complex, then $\beta(X)\tau i$ implies $X\tau_h i$.

Consequently if ${\mathcal F}$ is a family of maps with finitely dominated domains,
$\beta(X)\tau  {\mathcal F}$ implies $X\tau_h {\mathcal F}$.
\end{Lem}
\dokaz  
By ({\bf iii}) of Lemma \ref{ReductionIfTargetIsHomotopyDominated} we may assume that $K$ is in fact finite.

Let $f\colon A\to K$ be a map defined on the closed subset $A$ of $X$.
As $K$ is a finite CW complex, $f$ extends by \cite[Corollary $7.5.39$]{En$_1$}
to a map $F\colon\beta A\to K$. As $X$ is normal, $\beta A$ is the closure in
$\beta X$ of $A$. By assumption on $\beta X$, the map $i F$ extends to
a map $H\colon \beta X\to L$. Then $h=H\vert_X\colon X\to L$ is the desired
extension of $if$.
\hfill $\blacksquare$

\dokazp{Proposition \ref{quasi-finite}}
Follows immediately from Lemmas \ref{CechStoneCompLem} and \ref{ReverseOfCechStoneCompLem}.
\edokaz

\begin{Lem} \label{TechnicalLemma}
Let ${\mathcal F}$ be a family of maps between CW complexes and let $Y$ be a paracompact space.
Suppose that for every finite CW complex $M$ and every map $F\colon Y\to\Cone(M)$ there exists a factorization
$F=p\circ g$ through $g\colon Y\to Z$ where $Z$ is a paracompact space with $Z\tau_{cm}\FFF$. Then $Y\tau_{cm}\FFF$.
In particular, $Y\tau\FFF$ if $Y$ is compact.
\end{Lem}
\dokaz
Let $i\colon K\to L$ be a member of $\FFF$, and let $f\colon A\to K$ be a compact map defined on
the closed subset $A$ of $Y$. The image of $f$ is contained in a finite subcomplex $M$ of $K$.
Let $F\colon Y\to\Cone(M)$ denote an extension of the composite $A\xrightarrow{f}M\hookrightarrow\Cone(M)$.
Let $F=pg$ be the factorization guaranteed by the assumption. Then
\[ p\vert_{p^{-1}(M)}\colon p^{-1}(M)\to M\hookrightarrow K\xrightarrow{i}L \]
extends to a compact map $G\colon Z\to L$. The composite $Gg$ is the desired
extension of $if$.
\edokaz

\dokazp{\thmref{family_equivalent_to_map}}
Pick $t\in S$ and let $j_t\colon K_t\to\vee_sK_s$ and $q_t\colon\vee_sL_s\to L_t$ be the
obvious maps. Evidently $X\tau_hi$ implies $X\tau_hq_tij_t$. Therefore, $\setof{i}\lqs\FFF$ over paracompact spaces.

To show $\FFF\lqs\setof{i}$, let $X\tau_h\FFF$. In view of (3) and (4) of Proposition \ref{ResultsFromDy2}
it suffices to be able to extend any locally compact map $f\colon A\to\vee_sK_s$ defined
on an arbitrary closed subset $A$ of $X$ to a locally compact map defined on the whole space $X$.
We can extend the composite $A\xrightarrow{f}\vee_sK_s\to\vee_s\Cone(K_s)$
to a locally compact map $F\colon X\to\vee_s\Cone(K_s)$ by Corollary 2.10 of \cite{Dy2}. Let $*$ denote
the common base point of all the $K_s$. Set $A'=A\cup F^{-1}(*)$. By Corollary 2.11 of \cite{Dy2} we can
extend the restriction $F\vert_{A'}\colon A'\to\vee_sK_s$ to a locally compact map $f'\colon N\to\vee_sK_s$
for a closed neighborhood $N$ of $A'$. Define $N_s={f'}^{-1}(K_s)$ and $X_s=F^{-1}(\Cone (K_s))$. Note that
for distinct $t$ and $t'$, the intersection $X_t\cap X_{t'}$ equals $F^{-1}(*)$ where $*$ is the common
base point of all the $K_s$. By assumption and (4) of Proposition \ref{ResultsFromDy2}, the locally compact map
\[ N_s\xrightarrow{f'\vert{N_s}}K_s\xrightarrow{i_s}L_s	\]
extends to a locally compact map $g_s\colon X_s\cup N_s\to L_s$.
By construction, the $g_s$ define a function $g\colon X\to\vee_s L_s$ that extends $f'$.
Let $B$ be a closed neighborhood of $A'$ contained in the interior of $N$.
As $X\setminus B=\cup (X_s\setminus B)$ is the disjoint union of open sets, $g$ is continuous.
\edokaz

\dokazp{\thmref{EquivOverCompactaImpliesOverCompact}}
Let $X$ be a compact Hausdorff space with $X\tau {\mathcal F}$.
Let $M$ be a finite CW complex and $F\colon X\to\Cone(M)$ any map.
By Theorem \ref{FactorizationThmForFamily} the map $F$ factors as
$F=p\circ g$ through a surjective map $g\colon X\to Z$ where $Z$
is a compact metric space with $Z\tau\FFF$. By assumption this
implies also $Z\tau\GGG$. Lemma \ref{TechnicalLemma} implies $X\tau\GGG$.
\hfill $\blacksquare$

\dokazp{\thmref{TheQFThm}}
(1) Suppose that $\FFF$ is $\CCC$-quasi-finite and $X\in\CCC$. Evidently $X\tau_{cm}\FFF$
and hence $\beta(X)\tau\FFF$ by Proposition \ref{AEForCompactMapsAndCechStoneMap}.

(2) Suppose that $\FFF$ is range dominated by $\alpha$ and that $\CCC$ is $\alpha$-saturated.
By (1) we have to show that if $\FFF$ is weakly $\CCC$-quasi-finite then $\FFF$ is $\CCC$-quasi-finite.
To this end, suppose that for a member $i\colon K\to L$ of the family $\FFF$ the function $e=e_i$
of (2) of \defref{QFNotions} does not exist. Let $M_0$ be the smallest (finite) subcomplex of $L$
containing $i(M)$ and let $\{M_t\,\vert\,t\in T\}$ be all finite subcomplexes of $L$ that contain $M_0$. 
For each $t$ there are a space $X_t\in\CCC$ with $X_t\tau_h\FFF$, and a map $f_t\colon A_t\to M$
defined on the closed subset $A_t$ of $X_t$ such that there is no extension of the composite
$i\circ f_t$ to a map $X_t\to M_t$. Set $X=\oplus_{t\in T}X_t$, $A=\oplus_{t\in T}A_t$, and let
$f\colon A\to M$ be induced by $\{f_t\}_{t\in T}$. Obviously, $A$ is closed in $X$,
and, by assumption on $\FFF$ and $\CCC$, also $X\in\CCC$. Evidently $X\tau_h\FFF$, hence $X\tau_{cm}\FFF$,
by Proposition \ref{AEForCompactMapsAndCechStoneMap}. This means that for some $t$ there exists $g\colon X\to M_t$
with $g\vert_A=i\circ f$. But then $g_t=g\vert_{X_t}\colon X_t\to M_t$ has $g_t\vert_{A_t}=i\circ f_t$; contradiction.

(3) Suppose $\FFF$ is strongly $\CCC$-quasi-finite where $\CCC$ is rich.
Then $\FFF$ is $\CCC$-equivalent to a family $\GGG$ of maps between finite CW complexes.
We want to show that $X\tau_h\FFF$ implies $\beta(X)\tau\FFF$ to infer that $\FFF$ is weakly $\CCC$-quasi-finite.
To this end, let $X\in\CCC$ satisfy $X\tau_h\FFF$. As $\FFF\lqs\GGG$ over $\CCC$, this implies $X\tau_h\GGG$.
By \lemref{CechStoneCompLem} it follows that $\beta(X)\tau\GGG$. Because $\CCC$ is rich and $\GGG\lqs\FFF$ over $\CCC$,
\thmref{EquivOverCompactaImpliesOverCompact} yields $\beta(X)\tau\FFF$, as asserted.
\edokaz

\dokazp{\thmref{EquivOverCompAreEquiv}}
Suppose $\FFF$ and $\GGG$ are two strongly $\CCC$-quasi-finite families
of maps between CW complexes. This means that $\FFF$ and $\GGG$ are $\CCC$-equivalent
to families $\FFF_1$ and $\GGG_1$, respectively, of maps between finite CW complexes.

We want to show that if $\CCC$ is rich then $\FFF\iff\GGG$ over compact metric spaces
implies $\FFF\iff\GGG$ over $\CCC$. Let $\FFF\iff\GGG$ over compact metric spaces. Then
also $\FFF_1\iff\GGG_1$ over compact metric spaces. By \thmref{EquivOverCompactaImpliesOverCompact},
$\FFF_1\iff\GGG_1$ over the class of all compact Hausdorff spaces. 

Pick $X\in\CCC$ with $X\tau_h\FFF$. Then $\FFF\lqs\FFF_1$ over $\CCC$ implies $X\tau_h\FFF_1$.
\propref{quasi-finite} implies $\beta(X)\tau\FFF_1$ and since $\FFF_1\lqs\GGG_1$ over compact
Hausdorff spaces, also $\beta(X)\tau\GGG_1$. Another application of \propref{quasi-finite} yields
$X\tau_h\GGG_1$. Finally, $\GGG_1\lqs\GGG$ over $\CCC$ implies $X\tau_h\GGG$, as claimed.
\edokaz

\dokazp{\thmref{CharOfQFFamilies}}
(1) Assume that $\FFF$ is strongly $\CCC$-quasi-finite. Then $\FFF$ is $\CCC$-equivalent
to a family $\GGG$ of maps between finite CW complexes. By \thmref{EquivOverCompactaImpliesOverCompact},
the richness of $\CCC$ implies that the families $\FFF$ and $\GGG$ are equivalent over compact Hausdorff spaces.
The hypothesis is that $\FFF$ is $\CCC$-quasi-finite and hence weakly $\CCC$-quasi-finite by (1) of \thmref{TheQFThm}.
By definition, this means that $X\tau_h\FFF$ implies $\beta(X)\tau\FFF$ whenever $X\in\CCC$.
For the reverse implication, assume that $\beta(X)\tau\FFF$. Then $\beta(X)\tau\GGG$ because $\FFF\lqs\GGG$ over
compact Hausdorff spaces. This implies $X\tau_h\GGG$, by \propref{quasi-finite}. Therefore $X\tau_h\FFF$.
since $\GGG\lqs\FFF$ over $\CCC$.
\par (2) Assume now that $X\tau_h \FFF$ is equivalent to $\beta(X)\tau\FFF$ whenever $X\in\CCC$.
Let $\setof{e_i\,\vert\,i\in\FFF}$ be the collection of functions guaranteed
by (2) of \defref{EquFiniteOfFamiliesDef}. For each $i\colon K\to L$ that belongs to $\FFF$
and each finite subcomplex $M$ of $K$ let $j_{i,M}\colon M\to e_i(M)$ denote the restriction
of $i$ to $M$. Let $\GGG$ denote the collection of all maps $j_{i,M}$ where the index $i$
ranges over $\FFF$ and, for each $i\in\FFF$, the index $M$ ranges over all finite subcomplexes
of the domain of $i$. We want to show that $\FFF\iff\GGG$ over $\CCC$. The fact that $\FFF\lqs\GGG$ over $\CCC$
is the content of (2) of \defref{QFNotions}. To show that $\GGG\lqs\FFF$ over $\CCC$, assume that $X\in\CCC$
satisfies $X\tau_h\GGG$. By construction of $\GGG$ this implies $X\tau_{cm}\FFF$, and hence
$\beta(X)\tau\FFF$ by \propref{AEForCompactMapsAndCechStoneMap}. By assumption, this is
equivalent to $X\tau_h\FFF$.
\edokaz

\dokazp{\propref{ExistenceOfInvertibleImpliesWQF}}
Let $M$ be a finite CW complex, let $X\in\CCC$, and assume that $X\tau_h\FFF$.
Furthermore, let $F\colon X\to\Cone(M)$ be any map. We may assume that $\Cone(M)\subset I^\omega$,
and as $X\tau_h\FFF$, it follows that $F$ admits a lift $G\colon X\to p^{-1}(\Cone(M))$.
As $p^{-1}(\Cone(M))$ is compact, $p^{-1}(\Cone(M))\tau_{cm}\FFF$,
and \lemref{TechnicalLemma} (together with \propref{AEForCompactMapsAndCechStoneMap})
can be applied to deduce $\beta(X)\tau\FFF$.
\edokaz

\dokazp{\thmref{ExistenceOfInvertibleMapsForQFThm}}
Label all maps $Y\to I^\alpha$ where $Y\subset I^\alpha $ is closed and satisfies $Y\tau {\mathcal F}$
as $f_t\colon Y_t\to I^\alpha$, where $t$ belongs to the indexing set $T$. Let $Y$ be the \v Cech-Stone
compactification of $Y'=\bigoplus\limits_{t\in T} Y_t$. Let $\beta$ be the cardinality of $T$. Since $\FFF$
is $\CCC$-quasi-finite, $\FFF$ is quasi-finite with respect to $\beta$-saturation $\Bar\CCC^\beta$ of $\CCC$,
by \propref{saturation}. By \thmref{ExistenceOfSupComplexForFamilies}, $\FFF$ is also weakly $\Bar\CCC^\beta$-quasi-finite.
Therefore $Y\tau\FFF$, and \thmref{FactorizationThmForFamily} implies that the natural map $f\colon Y\to I^\alpha$
factors as $f=p\circ g$ where $g\colon Y\to E$ is a surjection onto a compact Hausdorff space $E$ with $E\tau\FFF$
and $w(E)\leq\alpha$. By construction of $Y$, the map $p$ is $\FFF$-invertible with respect to $\CCC$.
Indeed, let $X\in\CCC$ and $X\tau_h\FFF$. Any $u\colon X\to I^\alpha$ extends over $\beta(X)$, and
factors through some $f_t$ by \thmref{FactorizationThmForFamily}. Since all $f_t$ have lifts to $E$,
so has $u$.
\edokaz

\dokazp{\thmref{CharOfQFViaSepMetric}}
That (3) is equivalent to (2) if $\CCC$ is $\alpha$-saturated
follows from (2) of \thmref{ExistenceOfSupComplexForFamilies}.
\par That (1)$\implies$(3) if $\CCC$ is rich follows from (3) of \thmref{StronglyQFAreWQF}.
\par We show that (2)$\implies$(1). Let $\GGG$ be the collection of inclusions $j_M\colon M\to e(M)$
where $M$ ranges over all finite subcomplexes of $K$, such that $X\tau_h\GGG$ if $X\in\CCC$ satisfies $X\tau_h K$.
Suppose that $X\in\CCC$ satisfies $X\tau\GGG$. Theorem 2.9 of \cite{Dy2} says that $X\tau_{lc}K$ which,
according to Corollary 2.13 of \cite{Dy2}, amounts to $X\tau_h K$. Thus $K$ is $\CCC$-equivalent to $\GGG$ and $K$
is strongly $\CCC$-quasi-finite.
\edokaz

\dokazp{\corref{ModOfOriginalDefOfQF}}
Only (a)$\implies$(b) is non-trivial. 
Replace $K$ by a locally finite countable simplicial complex,
and let $X$ be a paracompact space. By (4) of \propref{ResultsFromDy2},
$X\tau_h K$ implies $X\tau K$. Given $f\colon X\to \Cone(M)$ where $M$ is
a finite CW complex, we apply Chigogidze's Factorization
Theorem \cite{Ch2} to obtain a factorization $F=p\circ G$ through
$G\colon X\to Y$ where $Y$ is a Polish space with $Y\tau K$.
Therefore $Y\tau_{cm}K$ and \ref{TechnicalLemma} says $X\tau_{cm}K$.
Thus $K$ is weakly quasi-finite with respect to paracompact spaces
and \thmref{CharOfQFViaSepMetric} implies that $K$ is quasi-finite
with respect to paracompact spaces.
\edokaz

\dokazp{\corref{QFAndCompactificationsOfSeparable}}
(a)$\implies$(b). Let $E$ be a compact metric space with $E\tau K$
and let $p\colon E\to I^\omega$ be an invertible map with respect to the class
of paracompact spaces $Y$ satisfying $Y\tau_h K$ (see \ref{ExistenceOfInvertibleMapsForQFThm}
and \ref{ModOfOriginalDefOfQF}). Given a compactification $\nu(X)$ of a separable metric space $X$
satisfying $X\tau_h K$ let $i\colon \nu(X)\to I^\omega$ be an embedding.
Lift $i\vert_X\colon X\to I^\omega$ to a map $j\colon X\to E$. The closure
of $j(X)$ in $E$ is the desired compactification of $X$.

\par (b)$\implies$(a). Suppose $X$ is a separable metric space satisfying $X\tau_h K$,
and let $f\colon X\to\Cone(M)$ be a map to the cone over a finite CW complex $M$.
Note that $f$ factors as $X\to \nu(X)\to\Cone(M)$ for some metric
compactification $\nu(X)$ of $X$. Therefore $f$ factors as $X\to\gamma(X)\to\Cone(M)$
for some metric compactification $\gamma(X)$ satisfying $\gamma(X)\tau K$.
Hence $X\tau_{cm}K$ by \lemref{TechnicalLemma}, and \thmref{CharOfQFViaSepMetric}
implies that $K$ is quasi-finite with respect to the class of separable metric spaces.
\edokaz

\dokazp{\thmref{WedgeOfQFIsQF}}
For the case of $\CCC$-quasi-finiteness suppose that the $K_s$ are $\CCC$-quasi-finite. Let
$\setof{e_s\,\vert\,s\in S}$ be the set of functions guaranteed by (2) of \defref{QFNotions}.
If $M$ is a finite subcomplex of $K$ then $M$ is contained in a finite wedge $\vee_{k=1}^nM_{s_k}$
of finite subcomplexes $M_{s_k}$ of $K_{s_k}$. Let $e(M)=\vee_{k=1}^ne_{s_k}(M_{s_k})$.
Let $X\in\CCC$ satisfy $X\tau_hK$. Then also $X\tau_hK_s$ for all $s$, and, in particular
$X\tau(M_{s_k}\to e_{s_k}(M_{s_k})$ for all $k$. By \thmref{family_equivalent_to_map}, it follows
that $X\tau(\vee_{k=1}^nM_{s_k}\to e(M))$ and therefore also $X\tau(M\to e(M))$. This shows
that $e$ satisfies (2) of \defref{QFNotions}, so $K$ is $\CCC$-quasi-finite.
\par For the case of strong $\CCC$-quasi-finiteness, let the $K_s$ be strongly $\CCC$-quasi-finite.
Then each $K_s$ is equivalent to a family, say $\FFF_s$, of maps between finite CW complexes.
By \propref{RealizationOfUnionOfFamilies}, the wedge $K=\bigvee_sK_s$ is equivalent to the union $\cup_s\FFF_s$.
Thus $K$ is also strongly $\CCC$-quasi-finite.
\par The remaining case is that of weak $\CCC$-quasi-finiteness. If the $K_s$ are weakly $\CCC$-quasi-finite,
then note that $\FFF=\setof{\id_s\colon K_s\to K_s}$ is a weakly $\CCC$-quasi-finite family. Indeed, if $X\in\CCC$
then $X\tau_h\FFF$ is equivalent to $X\tau_hK_s$ for all $s$, and this in turn implies $\beta(X)\tau K_s$ for all $s$
as the $K_s$ are $\CCC$-quasi-finite. By \thmref{family_equivalent_to_map}, the family $\FFF$ is equivalent
to $\setof{\id\colon K\to K}$ over the class of all paracompact spaces and therefore, by
\thmref{EquivOverCompactaImpliesOverCompact} also over the class of compact Hausdorff spaces.
Thus $\beta(X)\tau\FFF$ implies $\beta(X)\tau(\id\colon K\to K)$, as asserted.
\edokaz

\dokazp{\thmref{EveryQFIsWedgeOfCount}}
Let $e$ be the function defined on finite subcomplexes of $K$ that exists according to (2) of \defref{QFNotions}.
\par
By induction, we define a new function $E$ that also meets the requirements of (2) of \defref{QFNotions} and
has the additional property that $E(L)\subset E(L')$ provided $L\subset L'$. If $L$ has only one cell, we put $E(L)=e(L)$.
The general induction step is as follows. Add $E(M)$ to $e(L)$ for all proper $M\subset L$ and apply $e$ to that union.
\par In the next step we construct a family $K_M$ of countable subcomplexes of $K$ indexed by
all finite subcomplexes $M$ of $K$. Put $M_0=M$ and define $M_i$ inductively by $M_{i+1}=E(M_i)$.
Let $K_M=\bigcup\limits_{i=0}^\infty M_i$. Note that $M\subset L$ implies $K_M\subset K_L$.
Also, by Theorem 2.9 of \cite{Dy2}, $X\tau_h K$ implies $X\tau_h K_M$ for all $X\in\CCC$.
\par Consider $P=\bigvee\limits_{M\subset K}K_M$. \thmref{family_equivalent_to_map} implies that $K\lqs P$ over $\CCC$.
An application of Theorem 2.9 of \cite{Dy2} yields $P\leq K$ over the class of paracompact spaces.
Thus $P$ and $K$ are equivalent over $\CCC$. One can easily check that $P$ is $\CCC$-quasi-finite
using the function $e$.
\edokaz

%
%

\dokazp{\thmref{UnionThmForMaps}}
We may assume that $X_1\cup X_2=X$ and that all $K_1$, $K_2$, $L_1$, $L_2$ are simplicial complexes
equipped with the CW topology. For a simplicial complex $M$ let $|M|_w$ and $|M|_m$ denote the underlying
topological spaces with the CW (`weak') topology and the metric topology respectively. Any continuous map
$f\colon Y\to|M|_w$ is locally compact (see \cite{Dy2}). In addition, the following is true. If $f\colon Y\to|M|_m$
is a continuous map such that every $y\in Y$ has a neighborhood $U$ for which $f(U)$ is contained in a finite
subcomplex of $|M|_m$, then $f$ is also continuous as a function from $Y$ to $|M|_w$.

Suppose that $C$ is a closed subset of $X_1\cup X_2$ and let $f\colon C\to K_1\ast K_2$ be a map.
Note that $f$ defines two closed, disjoint subsets $C_1=f^{-1} (K_1)$,  $C_2=f^{-1} (K_2)$ of
$C$ and maps $f_1\colon C\setminus C_2\to K_1$, $f_2\colon C\setminus C_1\to K_2$, $\alpha\colon C\to [0,1]$ such that: 
\begin{enumerate}
\item $\alpha ^{-1} (0)=C_1$, $\alpha^{-1} (1)=C_2$, 
\item $f(x)=(1-\alpha (x))\cdot f_1(x)+\alpha (x)\cdot f_2(x)$ for all $x\in C$. 
\end{enumerate}

\par\noindent
Indeed, each point $x$ of a simplicial complex $M$ can be uniquely written as
$x=\sum\limits_{v\in M^{(0)}}\phi_v(x)\cdot v$ where $M^{(0)}$ is the set of vertices
of $M$ and $\{\phi_v(x)\}$ are barycentric coordinates of $x$. We define
$\alpha(x)=\sum\limits_{v\in K_2^{(0)}}\phi_v(f(x))$, and
\[ f_1(x)=\frac{1}{1-\alpha(x)}\sum\limits_{v\in K_1^{(0)}}\phi_v(f(x))\cdot v,\quad
	f_2(x)=\frac{1}{\alpha(x)}\sum\limits_{v\in K_2^{(0)}}\phi_v(f(x))\cdot v. \]
\par Since $X_1\setminus C_2\tau j_1$, the composite $j_1\circ f_1$ extends over $(C\cup X_1)\setminus C_2$. 
Consider a homotopy extension $g_1\colon U_1\to L_1$ of that map over a neighborhood $U_1$ of
$(C\cup X_1)\setminus C_2$ in $X\setminus C_2$. Since $C\setminus C_2$ is closed in $U_1$,
we may assume that $g_1$ is a genuine extension of $j_1\circ f_1\colon C\setminus C_2\to L_1$
(see Corollary 2.13 of \cite{Dy2}). Similarly, let $g_2\colon U_2\to L_2$ be
an extension of $j_2\circ f_2$ over a neighborhood $U_2$ of $(C\cup X_2)\setminus C_1$ 
in $X\setminus C_1$. Note that $X=U_1\cup U_2$. Let $\beta\colon X\to [0,1]$ be an extension
of $\alpha$ such that $\beta(X\setminus U_2)\subset\{0\}$ and $\beta(X\setminus U_1)\subset\{1\}$. 
Define  $f'\colon X\to L_1\ast L_2$ by
$$f'(x)=(1-\beta (x))\cdot g_1(x)+\beta (x)\cdot g_2(x)\text{ for all } x\in U_1\cap U_2,$$
$$f'(x)=g_1(x)\text{ for all } x\in U_1\setminus U_2,$$
and $$f'(x)=g_2(x)\text{ for all } x\in U_2\setminus U_1.$$ 
Note that $f'$ is a pointwise extension of $(j_1\ast j_2)\circ f$.  
To finish the proof it suffices to show that $f'\colon X\to |L_1\ast L_2|_m$
is continuous. Assuming this, the composite of $f'$ and a homotopy inverse
to the identity map $|L_1\ast L_2|_w\to |L_1\ast L_2|_m$ is an extension
of $f\colon C\to |L_1\ast L_2|_w$ up to homotopy. This implies the existence
of a genuine extension by Corollary 2.13 of \cite{Dy2}.
\par
To prove continuity of $f'\colon X\to |L_1\ast L_2|_m$ we need to show that the composites $\phi_v\circ f'$
are continuous for all vertices $v$ of $L_1\ast L_2$ (see Theorem 8 on page 301 in \cite {M-S}).
Without loss of generality we may assume that $v\in L_1$. Then, 
$$\phi_v(f'(x))=(1-\beta (x))\cdot \phi_v(g_1(x)) \text{ for all } x\in U_1$$
and $$\phi_v(f'(x))=0\text{ for all } x\in U_2\setminus U_1.$$ 
Clearly, the restriction $\phi_v\circ f'|_{U_1}$ is continuous. If $x_0\in U_2\setminus U_1$
is the limit of a sequence $\setof{x_n\vert n}$ contained in $U_1$ then the sequence $\setof{\beta(x_n)}$
converges to $1$ and $0\le\phi_v(g_1(x_n))\le 1$ for all $n$. Consequently, the sequence $\setof{\phi_v(f'(x_n))}$
converges to $0=\phi_v(f'(x_0))$ which finishes the proof.
\edokaz

\dokazp{\propref{RealizationOfUnionOfFamilies}} 
By assumption, $X \tau_h K_s\Leftrightarrow X \tau_h \mathcal F_s$ for each $s\in S$ and each $X\in\CCC$.
We are asserting $X\tau_h \bigvee_{s \in S}K_s \Leftrightarrow X \tau_h \bigcup_{s \in S}\mathcal F_s$.

Is $X \tau_h \bigvee_{s \in S}K_s$ then $X \tau_h K_s$ for all $s \in S$,
which means $X \tau_h \mathcal F_s$ and of course $X \tau_h \bigcup_{s \in S}\mathcal F_s$.

For the reverse implication, $X \tau_h \bigcup_{s \in S}\mathcal F_s$ implies
$X \tau_h \mathcal F_s$ for all $s \in S$ and therefore $X \tau_h\{\id\colon K_s\to K_s\}$.
Therefore $X \tau_h \bigvee_{s \in S}K_s$, by \thmref{family_equivalent_to_map}.
\edokaz

\dokazp{\thmref{JoinOfQF}}\label{ProofJoinOfQF}
Let $e_K$ be a function defined on the family of all finite subcomplexes of $K$ such that
$M\subset e_K(M)$ for all $M$ and $X\tau i_M$ for all Polish spaces $X$ satisfying $X\tau K$,
where $i_M$ is the inclusion $M\to e_K(M)$. Let $e_L$ be the analogous function for $L$.
Define $e$ for $K\ast L$ as follows: given a finite subcomplex $M$ of $K\ast L$
find finite subcomplexes $K_0$ of $K$ and $L_0$ of $L$ so that $M\subset K_0\ast L_0$.
Define $e(M)$ as $e_K(K_0)\ast e_L(L_0)$.
\par
Suppose $X$ is a Polish space and $X\tau K\ast L$. Express $X$ as $X_1\cup X_2$
so that $X_1\tau K$ and $X_2\tau L$ (see the main result of \cite{DrDy1}).
Notice $Y\tau K$ for every subset $Y$ of $X_1$ and $Y\tau L$ for every subset $Y$ of $X_2$.
Therefore (see \ref{UnionThmForMaps}) $X\tau i_M\ast i_P$ for every finite subcomplex $M$ of $K$
and every finite subcomplex $P$ of $L$. That readily implies $X\tau i_M$ for every finite
subcomplex $M$ of $K\ast L$, i.e. $K\ast L$ is quasi-finite.
\edokaz

\dokazp{\thmref{QFDimRealizedByCompacta}}
Consider a $K$-invertible map $p\colon X\to I^\omega$ to the Hilbert cube $I^\omega$
such that $X\tau K$ and $X$ is a compact metric space (such a map exists by \corref{InvExistCor}).
Suppose $X\tau L$ for some CW complex $L$. We need to show that $K\leq L$ over the class
of compact metric spaces (see \thmref{EquivOverCompAreEquiv}). Suppose that $Y\tau K$ where $Y$
is a compact metric space. We may assume that $Y\subset I^\omega$ in which case there is a lift
$g\colon Y\to X$ of the inclusion $Y\to I^\omega$. Thus $Y$ is homeomorphic to $g(Y)$, and hence
$Y\tau L$.
\edokaz

\dokazp{\propref{PropsOfQFDim}}
(1) is obvious as $X\tau_h K$ is equivalent to $\beta(X)\tau K$ for all
quasi-finite CW complexes by \thmref{CharOfQFViaSepMetric} and \thmref{CharOfQFFamilies}.
\par (2) Pick a quasi-finite complex so that $K=\dim_{QF}(X)$.
By \corref{InvExistCor} there exists a $K$-invertible map $p\colon E\to I^\omega$
from a compact metric space $E$ with $E\tau K$ to the Hilbert cube.
As in the proof of \thmref{QFDimRealizedByCompacta}, the equality $\dim_{QF}(E)=K$ holds.
Since $X$ embeds in $I^\omega$, a lift of that embedding is an embedding
of $X$ in $E$. Let $c(X)$ be the closure of $X$ in $E$. 
Since $\dim_{QF}(X)\leq \dim_{QF}(c(X))\leq \dim_{QF}(E)=K=\dim_{QF}(X)$,
we infer $\dim_{QF}(X)=\dim_{QF}(c(X))$.
\par (3) By (2) there exists a metric compactification $c(X)$ of $X$ of the same quasi-finite dimension as $X$.
By Lusin's Theorem there is a $G_\delta$-subset $X'$ of $Y$ containing $X$ so that $X'$
is homeomorphic to a subset of $c(X)$. Thus $\dim_{QF}(X)\leq \dim_{QF}(X') \leq  \dim_{QF}(c(X))$.
\edokaz

%
%

\section{Open problems}\label{op}

\begin{Problem}\label{SepMetricReverseOfCechStoneCompLem}
Let $X$ be a separable metric space and let $\FFF$ be a family of maps between CW complexes.
Find characterizations of families $\FFF$ for which $\beta(X)\tau\FFF$ implies $X\tau\FFF$.
\end{Problem}

\corref{ModOfOriginalDefOfQF} raises the following issue:

\begin{Problem}\label{AreQFForSeparableQFForParacompact} 
Suppose $K$ is a CW complex that is quasi-finite with respect to separable metric spaces.
Is $K$ quasi-finite (with respect to paracompact spaces)?
\end{Problem}

Chigogidze \cite{Ch2} proved the following

\begin{Thm}\label{Chigogidze}
Suppose $K$ is a countable CW complexes. If $f\colon X\to Y$ is a map of paracompact
spaces where $Y$ is a Polish space and $X$ satisfies $X_\tau K$, then $f$ factors as
$f=p\circ g$ where $g\colon X\to Z$ is surjective, $Z$ is a Polish space, and $Z\tau K$.
\end{Thm}

In view of \thmref{Chigogidze} the following problem is natural:

\begin{Problem}\label{Conjecture - generalized Chigogidze}
Suppose $\FFF$ is a countable family of maps between countable CW complexes.
If $f\colon X\to Y$ is a map of paracompact spaces such that $Y$ is Polish and $X\tau_h\FFF$,
does $f$ factor as $f=p\circ g$, where $g\colon X\to Z$ is surjective, $Z$ is Polish, and $Z\tau {\mathcal F}$?
\end{Problem}

A positive answer to \ref{Conjecture - generalized Chigogidze} would
imply a positive answer to \ref{AreQFForSeparableQFForParacompact}.

\begin{Problem}\label{AreQFEuivToWedgeOfFiniteQ} 
Is there a countable quasi-finite CW complex $K$ that is not equivalent (over the class
of paracompact spaces) to a wedge of finitely dominated CW complexes?
\end{Problem}

\begin{Problem}\label{AreQFEuivToCountableQ} 
Is there a quasi-finite (with respect to the class
of paracompact spaces) CW complex $K$ that is not equivalent (over the class
of paracompact spaces) to a countable quasi-finite CW complex?
\end{Problem}

\begin{Problem}\label{AreQFOverPolishQFOverParacompactQ} 
Is there a quasi-finite (with respect to the class of separable metric spaces)
CW complex $K$ that is not quasi-finite with respect to the class of paracompact spaces?
\end{Problem}

\begin{Problem}\label{JoinOfQFComplexesQ} 
If $K_1$ and $K_2$ are not necessarily countable quasi-finite CW complexes,
is their join $K_1\ast K_2$ (with the CW topology) quasi-finite?
\end{Problem}

\begin{Problem}\label{KGOneQuestion} 
Suppose $K(G,1)$ is quasi-finite. Is it equivalent over compact metric spaces to $S^1$?
Is it equivalent over separable metric spaces to $S^1$?
\end{Problem}

\begin{Problem}\label{SPQuestion} 
Suppose $K$ is a countable CW complex whose infinite symmetric product $SP(K)$
is quasi-finite. Is $SP(K)$ equivalent over compact metric spaces to $S^1$?
Is it equivalent over separable metric spaces to $S^1$?
\end{Problem}

\begin{Problem}\label{DoGeomRealizExistQ} 
Is there a family $\FFF$ of maps between finite CW complexes without a geometric realization?
\end{Problem}

\begin{Problem}\label{WhenIsTelescopeQFQ} 
Characterize those maps $f\colon K\to K$ whose
infinite mapping telescope of $K\to K\to K\to\ldots$ is quasi-finite.
\end{Problem}

\begin{Problem}\label{TelescopeAsGeomRealizQ} 
Characterize those maps $f\colon K\to K$ whose
infinite mapping telescope of $K\to K\to K\to\ldots$ is a geometric realization of $f$ over paracompact spaces.
\end{Problem}

\end{document}